\documentclass[preprint,number,12pt]{elsarticle}

\usepackage[all]{xy}
\RequirePackage{epsf}
\RequirePackage{amsmath}
\RequirePackage[dvips]{epsfig}
\RequirePackage{indentfirst}
\RequirePackage{fancyhdr}
\RequirePackage{color}
\RequirePackage{graphicx}
\RequirePackage{amssymb}

\usepackage{pdfscreen}
\usepackage{amsfonts} 
\usepackage{graphics}                 
\usepackage{hyperref}                 
\usepackage{anysize}
\usepackage{setspace}
\usepackage{amsthm}
\usepackage{algorithm}
\usepackage{algorithmic}

\newtheorem{theorem}{Theorem}
\newtheorem{lema}{Lemma}[section]
\newtheorem{prop}{Proposition}[section]
\newtheorem{cor}{Corollary}[section]
\theoremstyle{definition}

\newtheorem{remark}{Remark}[section]
\newtheorem{dfn}{Definition}[section]

\newenvironment{dem}{\subparagraph{Proof:}}{\qed}
\newenvironment{dem1}{\subparagraph{Proof of Theorem 2:}}{\qed}
\newenvironment{dem2}{\subparagraph{Proof of Theorem 3:}}{\qed}
\newenvironment{dem3}{\subparagraph{Proof of Theorem 4:}}{\qed}
\newcommand{\Ce}{\mathcal{C}}
\newcommand{\La}{\Lambda}    
\newcommand{\la}{\lambda}      
\newcommand{\ep}{\epsilon}
\newcommand{\s}{\sigma}
\newcommand{\be}{\beta}
\newcommand{\al}{\alpha}
\newcommand{\om}{\omega}
\newcommand{\barl}{\bar{\lambda}}

\newcommand{\barr}{\bar{r}}
\newcommand{\barx}{\bar{x}}
\newcommand{\bary}{\bar{y}}
\newcommand{\barq}{\bar{q}}

\newcommand{\U}{\mathcal{U}}
\newcommand{\SSS}{\mathcal{S}}

\newcommand{\M}{\mathbb{M}}

\newcommand{\R}{\mathbb{R}}
\newcommand{\Q}{\mathbb{Q}}
\newcommand{\Z}{\mathbb{Z}}

\newcommand{\A}{\mathbb{A}}
\newcommand{\OO}{\mathcal{O}}

\newcommand{\B}{\mathcal{B}}

\newcommand{\pe}{\mathfrak{p}}

\journal{Journal of Number Theory}

\begin{document}

\begin{frontmatter}



\title{Spinor norm for skew-hermitian forms over quaternion algebras}


\author[1]{Luis Arenas-Carmona}
\author[2]{Patricio Quiroz}

\address[1]{Universidad de Chile, Facultad de Ciencias, Casilla 653, Santiago, Chile}
\address[2]{Universidad de Chile, Facultad de Ciencias, Casilla 653, Santiago, Chile}

\begin{abstract}
We complete all local spinor norm computations for quaternionic skew-hermitian forms over the field of rational numbers. Examples of class number computations are provided.
\end{abstract}

\begin{keyword}

Skew-hermitian forms \sep spinor norm \sep quaternion algebras.
\end{keyword}

\end{frontmatter}



\section{Introduction}

Let $K$ be a number field and let $D$ be a quaternion algebra over $K$ with canonical involution $q \mapsto \barq$. Let $V$ be a rank-$n$ free $D$-module. Let $h:V\times V\rightarrow D$ be a \textit{skew-hermitian form}, i.e., $h$ is $D$-linear in the first variable and satisfies $h(x,y)=-\overline{h(y,x)}$. A $D$-linear map $\phi:V\rightarrow V$ preserving $h$ is called an isometry. We denote by $\U_K$ (resp. $\U_K^+$) the unitary group of $h$ (resp. the special unitary group of $h$), i.e., the group of isometries (resp. isometries with trivial reduced norm) of $h$. Skew-hermitian forms share many properties of quadratic forms. In fact, if $D\cong \M_2(K)$, skew-hermitian forms in a rank-$n$ free $D$-module are naturally in correspondence with quadratic forms in the $2n$-dimensional $K$-vector space $PV$, for any idempotent matrix $P$ of rank 1 in $D$ \cite[\S 3]{A-C07}. In this case, the unitary group of $h$ is isomorphic to the orthogonal group of the corresponding quadratic form. On the other hand, $\U_K=\U_K^+$ when $D$ is a division algebra \cite[\S 2.6]{K}.

As in the quadratic case, the problem of determining if two skew-hermitian lattices in the same space are isometric or not can be approached by the theory of genera and spinor genera. A skew-hermitian lattice or $\OO_D$-lattice in $V$, where $\OO_D$ is a maximal order in $D,$ is a lattice $\La$ in $V$ such that $\OO_D\La=\La.$ The special unitary group $\U_K^+$ acts naturally on the set of $\OO_D$-lattices and the $\U_K^+$-orbits are called classes (strict). This action can be extended to the adelization $\U_\A^+$ of $\U_K^+$ \cite[\S 2]{A-C03}. The $\U_\A^+$-orbits are called genera. A genus of $\OO_D$-lattices can be defined as a set of locally isometric lattices in the same space, since there is no Hasse principle for skew-hermitian spaces \cite{A-C07}. Between the class and the genus of a lattice lies its spinor genus. Two lattices $M$ and $\La$ are in the same spinor genus if, replacing each by an isometric lattice if needed, we can find, for each place $\pe$, local isometries $\s_\pe\in \U_{K_\pe}^+$ with trivial spinor norm satisfying $\s_\pe M_\pe=\La_\pe$. The spinor norm $\theta_\pe:\U_{K_\pe}^+\rightarrow K_\pe^*/K_\pe^{*2}$, or more generally, $\theta_L:\U_L^+\rightarrow L^*/L^{*2}\cong H^1(L,F)$ for any field $L$ containing $K$, is the coboundary map defined from the universal cover $F_{\bar{L}}\hookrightarrow\U_{\bar{L}}^+\twoheadrightarrow\U_{\bar{L}}^+$ \cite[\S 2]{A-C03}. This concept is important because both the class and the spinor genus of a given lattice coincide whenever $\U_{K_\pe}^+$ is non-compact for some arquimedian place $\pe.$ The set of classes contained in the genus of a lattice $\La$ is in one-to-one correspondence with the set of double cosets $$\U_K^+\setminus \U_\A^+/\U_\A^+(\La),$$ 
where $\U_\A^+(\La)$ is the adelic stabilizer of $\La$. The cardinality of $\U_K^+\setminus \U_\A^+/\U_\A^+(\La)$ is called the \textit{class number} of  $\La$ with respect to $\U_K^+$. This quantity is difficult to compute in general. An easier problem is determine the number of spinor genera in a genus. In fact, this number is equal to the order of the finite abelian group $$\Theta_\A(\U_\A^+)/\Big(\theta(\U_K^+)\Theta_\A\big(\U_\A^+(\La)\big)\Big),$$ where $\Theta_\A$ is the adelic spinor norm \cite[\S 2]{A-C03}.
Moreover, if $P:J_K\rightarrow J_K/J_K^2$ is the natural projection, where $J_K$ is the idelic group of $K$, and $$H_\A(\La):=P^{-1}\Big(\Theta_\A\big(\U_\A^+(\La)\big)\Big),$$ we have the following group isomorphism \cite[\S 2]{A-C03}: $$\Theta_\A(\U_\A^+)/\Big(\theta(\U_K^+)\Theta_\A\big(\U_\A^+(\La)\big)\Big)\cong J_K/K^*H_\A(\La).$$ To compute the group on the right, we need to know the image of the local spinor norm $\theta_\pe:\U_{K_\pe}^+(\La)\rightarrow K_\pe^*/K_\pe^{*2}$ at each place $\pe$ of the number field $K$. This is why we are interested in local spinor norm computations. Full computations exist for symmetric integral bilinear forms. Non-dyadic cases can be found in \cite{K1} and dyadic cases in \cite{Beli}. For this reason we assume, from now on, that the quaternion algebra $D$ is a division algebra. Remember that in this case, we have $\U_K=\U_K^+.$ For skew-hermitian forms, non-dyadic places have been completely studied by B\"oge in \cite{B}. The dyadic case was studied by Arenas-Carmona in \cite{A-C04} and \cite{A-C10}, not completing all the cases, which we carry out here when $K_\pe=\Q_2$. From now on $k=K_\pe$ denotes a dyadic local field of characteristic 0. 

If $D$ is a division algebra over $k$ we can define an absolute value $|\cdot|:D\rightarrow \R_{\geq 0}$ by $|q|=|Nq|_k$, where $N$ is the reduced norm and $|\cdot|_k$ is the absolute value of $k.$ The valuation on $D$ induced by $|\cdot|$ is denoted by $\nu$. Any skew-hermitian lattice $\Lambda$ has a decomposition of the type \begin{equation}\Lambda = \Lambda_1\bot \cdots \bot \Lambda_t, \label{dec}\end{equation} where each lattice $\Lambda_r$ has rank 1 or 2, and the scales satisfy $\textbf{s}(\Lambda_{r+1})\subset \textbf{s}(\Lambda_r)$ \cite[\S 5]{A-C04}. This is the skew-hermitian analogue to the Jordan decomposition for bilinear lattices in \cite[\S 91]{O}. If some $\Lambda_m$ in the decomposition of $\La$ has rank 1, then $\Lambda_m=\OO_Ds_m$ and $h(s_m,s_m)=a_m$. Define $A\subset k^*/k^{*2}$ by $A=\{N(a_m)k^{*2}|\textup{ rank}(\Lambda_m)=1\}$. Following \cite{B}, we define $H(\Lambda)$ by the relation $H(\Lambda)/k^{*2}=\theta(\U_k^+(\Lambda)),$ where $\U_k^+(\Lambda)$ is the stabilizer of $\Lambda$ in $\U_k^+$, and $\theta:\U_k^+(\Lambda)\rightarrow k^*/k^{*2}$ denotes the spinor norm. Note that $k^{*2}\subset H(\Lambda)\subset k^*.$ By abuse of language we say that $H(\Lambda)$ is the image of the spinor norm. It is clear that if $\La=\La_1\bot\La_2$, then $H(\La_1),H(\La_2)\subset H(\La)$. In addition, for binary indecomposable lattices, we know that $H(\La)=k^*$ \cite{A-C10}. The lattices $\La$ for which the value of $H(\La)$ remains unknown to date are:\\

\noindent\textbf{Case I:} $\La=\langle a_1\rangle\bot\cdots\bot\langle a_n\rangle$, where $A=\{-uk^{*2}\}$ and the minimal difference between the valuation of the scales of two consecutive components satisfy $0<\min\{\nu(a_{i+1})-\nu(a_i)\}\leq \nu(16).$ Here, $u\in\OO_k^*$ denotes an arbitrary unit of non-minimal quadratic defect \cite[\S 63]{O}.\\
\textbf{Case II:} $\La=\langle a_1\rangle\bot\cdots\bot\langle a_n\rangle$, where $A=\{ \pi k^{*2}\}$ and the minimal difference between the valuation of the scales of two consecutive components satisfy $\nu(4)\leq\min\{\nu(a_{i+1})-\nu(a_i)\}\leq \nu(16).$ Here, $\pi$ denotes a prime in $k$.\\

In this article, we compute $H(\La)$ for cases I and II above, when $k=\Q_2$. Concretely, we have the following result: 

\begin{theorem}

The following table contains all local spinor norm computations when the base field is $\Q_2:$

\begin{table}[h]
\begin{center}
    $$\begin{array}{llllllp{5cm}}
    
    s & |A| & A & \mu & H(\Lambda) & \text{Reference}\\ \hline
     - & >1 & - & - & \Q_2^* & \text{Corollary } \ref{kest2}\\
    0 & 1 & \text{-}\Delta\Q_2^{*2} & - & \Z_2^*\Q_2^{*2} & \text{Table 2 in \cite{A-C04}}\\
    0 & 1 & \text{-}u\Q_2^{*2} & 0\leq\mu<\nu(8) & \Q_2^*& \text{Prop. } \ref{propcaso1k}+\text{Table 1 in\cite{A-C04}}\\
     0 & 1 & \text{-}u\Q_2^{*2}  & \mu\geq\nu(8) & N\big(\Q_2(a_m)^*\big)& \text{Prop. } \ref{propcaso1}+\text{Table 2 in\cite{A-C04}}\\
        0 & 1 & \pi\Q_2^{*2} & 0\leq\mu\leq\nu(16) & \Q_2^*& \text{Prop. } \ref{propcaso2}+ \text{Prop.} \ref{nu4}+\text{Tables 1,2 in\cite{A-C04}}\\
    0 & 1 & \pi\Q_2^{*2} & \mu>\nu(16) & N\big(\Q_2(a_m)^*\big)& \text{Table 2 in \cite{A-C04}}\\
    \neq 0 & - & - & - & \Q_2^*& \text{Theorem 2 in \cite{A-C10}}\\
    \end{array}$$
\end{center}
\caption{Spinor images for arbitrary lattices over $\Q_2$.}
\label{completa}
\end{table}

Here, $s$ denotes the number of indecomposable components of rank 2 in the decomposition (\ref{dec}) of $\La$, $\mu=\mu(\La)$ denotes the minimal difference between the valuation of the scales of two consecutive components of rank 1, and $\Delta\in \OO_k^*$ is a unit of minimal quadratic defect \cite[\S 63]{O}. Furthermore, $A$, $\pi$ and $u$ are as in the previous discussion. A dash means irrelevant information.
\label{teo1}
\end{theorem}

Our (computer assisted) proof of Theorem \ref{teo1} is based on the following scheme:

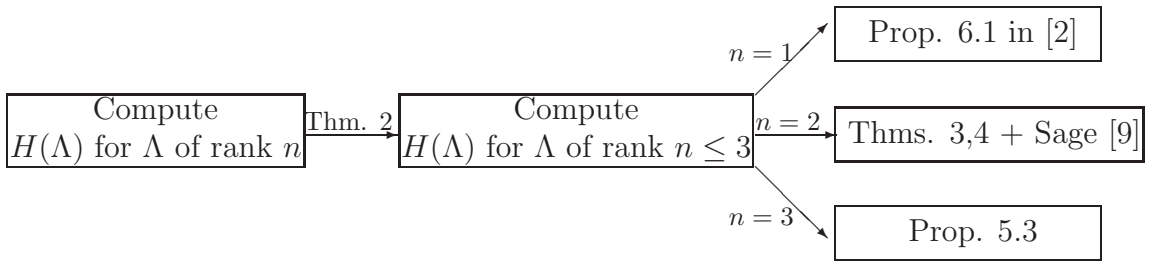
\begin{picture}(380,80)
\put(0,0){\framebox(111,27){\shortstack{Compute\\$H(\Lambda)$ for $\Lambda$ of rank $n$ }}} 

\put(111,12){\vector(1,0){35}}
\put(111,14){{\footnotesize Thm. 2}}
\put(147,0){\framebox(132,27){\shortstack{Compute\\$H(\Lambda)$ for $\Lambda$ of rank $n\leq 3$ }}}
\put(280,12){\vector(1,0){30}}
\put(280,14){{\footnotesize $n=2$}}
\put(310,2){\framebox(115,20){ Thms. 3,4 + Sage [9]}}
\put(280,27){\vector(1,1){27}}
\put(280,0){\vector(1,-1){27}}
\put(310, 40){\framebox(99,20){ Prop. 6.1 in [2]}}
\put(310, -35){\framebox(99,20){ Prop. 5.3}}
\put(270,40){{\footnotesize $n=1$}}
\put(270,-22){{\footnotesize $n=3$}}
\end{picture}

\vspace*{1.5cm}


 

The following result is useful to reduce the study of $H(\Lambda)$ to the case of low rank $\Lambda$ for arbitrary local fields. 

\begin{theorem}
Let $\Lambda=\langle a_1\rangle\bot\cdots\bot\langle a_n\rangle=\OO_Ds_1\bot\cdots\bot\OO_Ds_n$ be a skew-hermitian lattice and let $\mu=\mu(\La)$ be as above. Assume $\;\mu>\nu(4)$ and $N(a_2),...,N(a_n)\in N(a_1) k^{*2}$. Let $(s;\s)\in\B(\La)$, i.e., $s=(1-r)s_m-s_0$, where $s_0=\lambda_{m+1}s_{m+1}+\cdots+\lambda_ns_n\in \OO_Ds_{m+1}\bot\cdots\bot\OO_Ds_n$, $\s=a_m(1-\barr)$ and $|1-\barr|\geq |2|$. If $|\lambda_{m+t}|\geq |\la_{m+t+l}|$, for some $t\in\{1,...,n-m\}$ and for all $l\in\{1,...,n-m-t\}$, then there exists $\Lambda'=\langle b_1\rangle\bot\cdots\bot\langle b_{t+1}\rangle\subset\La$ satisfying the following conditions:
\begin{enumerate}
\item $(s;\s)\in\U_k^+(\La').$
\item $\mu(\La')\geq\mu(\La)$.
\item $N(b_i) \in N(a_1)k^{*2}$, for all $i=1,...,t+1$.
\end{enumerate}  
\label{red}
\end{theorem}

The following theorems help us to develop an algorithm to compute $H(\Lambda)$ for lattices of rank 2 for unramified dyadic fields. Using these results, together with Theorem \ref{red}, we compute $H(\La)$ for the unknown cases when the base field is $\Q_2$. Remember that, if $\La=\La_1\bot\La_2$, then $H(\La_1),H(\La_2)\subset H(\La)$. So, if we want to prove that $H(\La)=k^*$, is enough to do it for binary lattices.

\begin{theorem}
Let $\Lambda=\langle a_1\rangle\bot\langle a_2\rangle$ be a skew-hermitian lattice such that $|2a_1|\geq |a_2|$ and $N(a_2)\in N(a_1)k^{*2}$. The following statements are equivalent:
\begin{enumerate}
\item $H(\Lambda)=k^*$.
\item There exists $(s;\s)\in\B(\La)$ such that $N\s\notin N\big(k(a_1)^*\big).$
\item There exists $r\in\OO_D$ such that: $\left(\frac{N(1-r)\,,\,-Na_1}{\pe}\right)=-1$, $NzNa_1\in k^{*2}$ and $NzN(\pi^ta_1)^{-1}$ $\in \OO_k$, where $z=a_1-ra_1\barr$ and $\nu(\pi^t)=\mu(\La).$   
\end{enumerate}
\label{eq}
\end{theorem}

Let $e=\nu(2)/2$ be the ramification index of $k/\Q_2$ and remember that $t$ in the previous theorem satisfies $\nu(\pi^t)=\mu(\La)$. We call the conditions for $r$ in the statement \textit{3} of the Theorem \ref{eq}, the \textit{$k$-star conditions}.

\begin{theorem}
There exists $r\in \OO_D$ satisfying the $k$-star conditions if and only if there exists $\al\in \mathcal{S}\oplus\mathcal{S}\om\oplus\mathcal{S} i\oplus\mathcal{S} i\om\subset\OO_D$ satisfying them, where $\mathcal{S}$ is any finite set of representatives of $\OO_k/\pi^u\OO_k$, with $u=t+6e$.
\label{ssigen}
\end{theorem}

\section{Arithmetic of $D=\left(\frac{\pi,\Delta}{k}\right)$}\label{ad}

From here on, we work with the unique quaternion division algebra, up to isomorphism, $D=\left(\frac{\pi,\Delta}{k}\right).$ A base for $D$ is denoted by $\{1,i,j,ij\}$, where $$i^2=\pi,\;j^2=\Delta,\;ij=-ji.$$

Let $\OO_D$ be the unique maximal order in $D.$ The ring $\OO_D$ is the set of integral elements in $D$ \cite[\S 2]{V}, i.e., $\OO_D=\{q\in D \;|\; |q|\leq 1\}$. Moreover, we have $\OO_D=\OO_{k(j)}\oplus i\OO_{k(j)}$. Is not difficult to prove that $\OO_{k(j)}=\OO_k[\om],$ where $\om=\frac{j+1}{2}.$ Hence, we can write $\OO_D=\OO_k\oplus\OO_k\cdot\om\oplus\OO_k\cdot i\oplus\OO_k\cdot i\om,$ i.e., $\{1,\om,i,i\om\}$ is a base for $\OO_D.$ In general, if we denote by $\widetilde{\om}$ and $\widetilde{i}$ the classes, modulo $i^t$, of $\om$ and $i$ respectively, we have: {\footnotesize $$\OO_D/i^t\OO_D\cong \begin{cases} (\OO_k/\pi^s\OO_k) \oplus (\OO_k/\pi^s\OO_k)\widetilde{\om} \oplus (\OO_k/\pi^s\OO_k) \widetilde{i} \oplus (\OO_k/\pi^s\OO_k) \widetilde{i}\widetilde{\om} & , \text{ if }t=2s \textup{ even. }\\ (\OO_k/\pi^{s+1}\OO_k) \oplus (\OO_k/\pi^{s+1}\OO_k)\widetilde{\om} \oplus (\OO_k/\pi^s\OO_k) \widetilde{i} \oplus (\OO_k/\pi^s\OO_k) \widetilde{i}\widetilde{\om} & , \text{ if }t=2s+1 \textup{ odd. } \end{cases}$$}

\begin{remark}
If $\al=a+b\om+ci+di\om\in\OO_D$, then the trace $T$ and the reduced norm $N$ of $D$ satisfy $T(\al)=2a+b$ and $N(\al)=a^2+ab-\delta b^2-\pi(c^2+cd-\delta d^2)$, where $\delta\in \OO_k^*$ satisfies $\Delta=1+4\delta.$
\end{remark}

We finish this section with a key result that we use frequently in this article to prove that certain norms of quaternions are in the same square class.

\begin{lema}
If $\al\in\OO_D$ satisfies $|\al|<1$, then $N(1+4\al)$ is a square.
\label{LTCL}
\end{lema}

\begin{dem}
$N(1+4\al)=1+4T(\al)+16N(\al)=1+4(T(\al)+4N(\al))$. The condition $|\al|<1$ implies $\pi|T(\al)$. Hence, the result follows from the Local Square Theorem \cite[\S 63]{O}.
\end{dem}\\

\section{Generators of $\U_k^+(\Lambda)$ and their spinor norm}\label{secgen}

Let $(V,h)$ be a skew-hermitian space and let $s\in V,\;\s\in D^*$ be such that $\s-\bar{\s}=h(s,s).$  Define the map $(s;\s):V\rightarrow V$ by $$(s;\s)(x)=x-h(x,s)\s^{-1}s.$$ Following \cite{A-C04} we call such maps \textit{simple rotations} with axis of rotation $s$. Simple rotations generate $\U_k^+$ and satisfy $$\theta[(s;\s)]=N(\s)k^{*2} \;\text{ \cite{S-H}},$$ where $N:D^*\rightarrow k^*$ is the reduced norm. Note that $\s\in k(a)$, where $a=h(s,s).$\\

In \cite[\S 6]{A-C04} the following lemmas are proved. The first one tell us how to construct simple rotations and is used in the second to obtain a set of generators for $\U_k^+(\Lambda)$.

\begin{lema}\cite[Lemma 6.3]{A-C04}
Let $(V,h)$ be a skew-hermitian space, and let $t,u\in V$ be such that $h(u,u)=h(t,t)=a.$ Define $r$ and $t_0$ by $u=rt+t_0,$ where $t_0\in t^\bot.$ Let $s=t-u$ and $\s=h(t,s).$ Then the following identities hold:

$$\s=a(1-\barr), \hspace{0.3cm} h(t_0,t_0)=a-ra\barr,\hspace{0.3cm} \s-\bar{\s}=h(s,s).$$ 

In particular $(s;\s)$ is a well-defined simple rotation satisfying $(s;\s)(t)=u$.
\label{gen0}
\end{lema}

If $\phi\in\U_k^+$, we have $h(\phi(t),\phi(t))=h(t,t)$, and hence, there exists a simple rotation $(s;\s)$ such that $(s;\s)(t)=\phi(t).$ This fact can be used to prove the following result by an induction process.

\begin{lema}\cite[Lemma 6.7]{A-C04}
Let $\Lambda$ be as in (\ref{dec}). Assume that $\La_m=\OO_Ds_m$, with $a_m=h(s_m,s_m)$, for $m=1,...,n$. Assume also $|2a_m|\geq |a_l|$ for $m<l$. Then the unitary group $\U_k^+(\Lambda)$ of the lattice is generated by elements of the following types:\\A) Simple rotations with axis $s_m, \,$ for some $m=1,...,n.$ \\ B) Simple rotations of the form $(s;\s)$, where we have $\s=a_m(1-\barr),$ and $s=(1-r)s_m-s_0,$ for some $ s_0\in\OO_Ds_{m+1}\bot\cdots\bot\OO_Ds_n.$ Furthermore, we can assume $1-r\notin (2i).$                                                                                                                                                                                                                                                                                                                                                                                                                                                                                                                                                                                                                                                                                                                                                                                                                                                                             
\label{gen}
\end{lema}

\begin{dfn}
Let $\La$ be as in (\ref{dec}). We denote by $\mathcal{A}(\La)$ the set of simple rotations of type (A) and by $\B(\La)$ the set of simple rotations of type (B) such that $|1-r|\geq |2|$, i.e., such that $1-r\notin (2i)$. It follows that, $\U_k^+(\La)$ is generated by $\mathcal{A}(\La)\cup \B(\La).$
\label{defAB}
\end{dfn}

If $\Lambda=\langle a_1\rangle \bot ...\bot\langle a_n\rangle$ is a skew-hermitian lattice, then $[k^*:N\big(k(a_1)^*\big)]=2$ \cite[\S 63]{O} and $N\big(k(a_1)^*\big)\subset H(\Lambda)$ \cite[Proposition 6.1]{A-C04}. A direct consequence of these facts is the following result:

\begin{prop}
Let $\Lambda=\langle a_1\rangle\bot ...\bot\langle a_n\rangle$ be a skew-hermitian lattice. Then $H(\Lambda)=N\big(k(a_1)^*\big)$ or $H(\La)=k^*.$
\label{konok}
\end{prop}

\begin{cor}
Let $\La$ be as in the proposition. If there exists $b\in \OO_D$ with $N(b)\notin N(k(a_1)^*)$ such that $\La=\langle b\rangle\bot \La'$, then $H(\La)=k^*.$
\label{kest2}                                                                                                                                                                                                                                                                                                                                                                                                                                                                                                                                                                                                                                                                                                                                                                                                                                                                                                                                                                                                                                                                                                                                                                                                                                                                                                                                                                                                                                                                                                                                                                                                                                                                                                                                                                                                                                                                                                                                                                                                                                                                                                                                                                                                                                                                                                                                                                                                                                                                                                                                                                                                                                                                                                                                                                                                                                                                                                                                                                                                                                                                                                                                                                                                                                                                                                                                                                                                                                                                                                                                                                                                                                                                                                                                                                                                                                                                                                                                                  
\end{cor}

\begin{cor}                                                                                                                                                                                                                                                                                                                                                                                                                                                                                                                                                                                                                                                                                                                                                                                                                                                                                                                                                                                                                                                                                                                                                                                                                                                                                                                                                                                                                                                                                                                                                                                                                                                                                                                                                                                                                                                                                                                                                                                                                                                                                                                                                                                                                                                                                                                                                                                                                                                                                                                                                                                                                                                                                                                                                                                                                                                                                                                                                                                                                                                                                                                                                                                                                                                                                                                                                                                                                                                                                                                                                                                                                                                                                                                                                                                                                                                                                                                                                                                                                                                                                                 Let $\La$ be as in the proposition. Then $H(\Lambda)=k^*$ if and only if, there exists $\phi\in\Ce(\Lambda)$ such that $\theta(\phi)\notin N\big(k(a_1)^*\big)/k^{*2}$, where $\Ce(\La)$ is a set of generators for $\U_k^+(\La).$                                                                                                                                                                                                                                                                                                                                                                                                                                                                                                                                                                                                                                                                                                                                                                                                                                                                                                                                                                                                                                                                                                                                                                                                                                                                                                                                                                                                                                                                                                                                                                                                                                                                                                                                                                                                                                                                                                                                                                                                                                                                                                                                                                                                                                                                                                                                                                                                                                                                                                                                                                                                                                                                                                                                                                                                                                                                                                                                                                                                                                                                                                                                                                                                                                                                                                                                                                                                                                                                                                                                                                                                                                                                                                                                                                                                                                                                                                                                                   \label{kest}                                                                                                                                                                                                                                                                                                                                                                                                                                                                                                                                                                                                                                                                                                                                                                                                                                                                                                                                                                                                                                                                                                                                                                                                                                                                                                                                                                                                                                                                                                                                                                                                                                                                                                                                                                                                                                                                                                                                                                                                                                                                                                                                                                                                                                                                                                                                                                                                                                                                                                                                                                                               \end{cor}                                                                                                                                                                                                                                                                                                                                                                                                                                                                                                                                                                                                                                                                                                                                                                                                                                                                                                                                                                                                                                                                                                                                                                                                                                                                                                                                                                                                                                                                                                                                                                                                                                                                                                                                                                                                                                                                                                                                                                                                                                                                                                                                                                                                                                                                                                                                                                                                                                                                                                                                                                                                                                                                                                                                                                                                                                                                                                                                                                                                                                                                                                                                                                                                                                                                                                                                                                                                                                                                                                                                                                                                                                                                                                                                                                                                                                                                                                                                                                                                                                                                                                                                                                                                                                                                                                                                                                                                                                                                                                                                                                                                                                                                                                                                                                                                                                                                                                                                                                                                                                                                                                                                                                                                                                                                                                                                                                                                                                                                                                                                                                                                                                                                                                                                                                                                                                                                                                                                                                                                                                                                                       \begin{remark}
Simple rotations $(s_m;\s)\in\mathcal{A}(\La)$ have spinor norm $\theta[(s;\s)]=N(\s)k^{*2}\in N\big(k(a_m)^*\big)/k^{*2}$, since $\s\in k(a_m).$ Hence, if $|2a_m|\geq |a_l|$ for $m<l,$ in the previous corollary we can set $\Ce(\La)=\B(\La)$. 
\label{obscorkest}
\end{remark}                                                                                                                                                                                                                                                                                                                                                                                                                                                                                                                                                                                                                                                                                                                                                                                                                                                                                                                                                                                                                                                                                                                                                                                                                                                                                                                                                                                                                                                                                                                                                                                                                                                                                                                                                                                                                                                                                                                                                                                                                                                                                                                                                                                                                                                                                                                                                                                                                                                                                                                                                                                                                                                                                                                                                                                                                                                                                                                                                                                                                                                                                                                                                                                                                                                                                                                                                                                                                                                                                                                                                                                                                                                                                                                                                                                                                                                                                                                                                                                                                                                                                                                                                                                                                                                    

The following result is key to prove Theorem \ref{ssigen}. Note that, $H(\La)$ depends on the existence of simple rotations with specific spinor norm (see Corollary \ref{kest}).

\begin{lema}
Let $\La=\langle a_1\rangle\bot ...\bot\langle a_n\rangle=\OO_Ds_1\bot\cdots\bot\OO_Ds_n$ be a skew-hermitian lattice such that $|2a_m|\geq |a_l|$ for $m<l$. Take $(s;\s)\in\B(\La)$, i.e., $s=(1-r)s_m-s_0$, where $$s_0=\lambda_{m+1}s_{m+1}+\cdots+\lambda_ns_n\in \OO_Ds_{m+1}\bot\cdots\bot\OO_Ds_n,\;\s=a_m(1-\barr) \text{ and } |1-r|\geq|2|.$$ If any of the following conditions is satisfied:
\begin{enumerate}
\item $|1-r|>|2|$ and $|\la_{m+1}|<1$, for $\mu(\La)\geq\nu(8),$ and the extension $k/\Q_2$ is unramified, 
\item $|1-r|=|2|$ and $|\la_{m+1}|\leq |2|,$ for $\mu(\La)\geq\nu(4\pi),$
\item $|1-r|>|2|$ or $|\la_{m+1}|<1$, for $\mu(\La)\geq\nu(16),$
\item $m=1$ and $|\la_2|\leq |4|,$
\end{enumerate} 
then $\theta[(s;\s)]\in N\big(k(a_m)^*\big)/k^{*2}.$                                                                                                                                                                                                                                                                                                                                                                                                                                                                                                                                                                                                                                                                                                                                                                                                                                                                                                                                                                                                                                                                                                                                                                                                                                                                                                                                                                                                                                                                                                                                                                                                                                                                                                                                                                                                                                                                                                                                                                                                                                                                                                                                   
\label{v16}
\end{lema}

\begin{dem} It suffices to prove that if $a=h(s,s)$, then $N(a)\in N(a_m)k^{*2}$, since $\s\in k(a).$ In fact, we have $s=(1-r)s_m-s_0$, so that $a=(1-r)a_m(1-\barr)+a_0$, where $a_0=h(s_0,s_0).$ It follows that 
\begin{equation}
N(a)=N(a_m)N(1-r)^2N\big(1+(1-r)^{-1}a_0(1-\barr)^{-1}a_m^{-1}\big).
\label{ecnorma}
\end{equation} 
Now, 
 $a_0=\la_{m+1}a_{m+1}\overline{\la_{m+1}}+\cdots+\la_na_n\overline{\la_n}$ and $|(1-r)^{-1}a_0(1-\barr)^{-1}a_m^{-1}|=|a_0||a_m^{-1}|/|1-r|^2<|4|$ if any of the conditions above is satisfied. This implies that the last norm in (\ref{ecnorma}) is a square in virtue of Lemma \ref{LTCL}. 
\end{dem}

\section{Proof of Theorems 2,3 and 4}
\begin{dem1} Set $\Lambda'=\OO_Ds_m\bot\cdots\bot\OO_Ds_{m+t-1}\bot\OO_Ds_0'=$ $\langle b_1\rangle\bot\cdots\bot\langle b_{t+1}\rangle$, where $s_0'=\la_{m+t}s_{m+t}+\cdots+\la_ns_n$, $b_i=h(s_{m+i-1},s_{m+i-1})=a_{m+i-1},$ for $i=1,...,t$ and $b_{t+1}=h(s_0',s_0')$. It is clear that $\La'\subset \La$. To prove the condition \textit{1} in the theorem is satisfied, we note that $s_0=s_0'-(\la_{m+1}s_{m+1}+\cdots+\la_{m+t-1}s_{m+t-1})\in\La'$. We compute, $$(s;\s)(s_m)=rs_m+s_0\in\Lambda',$$ $$(s;\s)(s_i)=s_i-h(s_i,s)\s^{-1}s=s_i+h(s_i,s_0)\s^{-1}s, \text{ for } i=m+1,...,m+t-1,$$ and $$(s;\s)(s_0')=s_0'-h(s_0',s)\s^{-1}s=s_0'+h(s_0',s_0)\s^{-1}s.$$ Hence, $(s;\s)(s_i),(s;\s)(s_0')\in\La'$ if $h(s_i,s_0)\s^{-1},h(s_0',s_0)\s^{-1}\in\OO_D$. The latter holds since $|\s|=|a_m(1-\barr)|\geq |2a_m|\geq \mathbb{S}(s_0),\mathbb{S}(s_0')$, where $\mathbb{S}(v)=\max_{w\in\La}|h(v,w)|$ is the height of a vector $v\in\La$ \cite[\S 5]{A-C04}. Then, $(s;\s)\in\U_k^+(\La')$. On the other hand, as $$b_{t+1}=h(s_0',s_0')=\sum_{u=m+t}^n \lambda_ua_u\overline{\la_u},$$ we have $|b_{t+1}|=|a_{m+t}||\la_{m+t}|^2$, since $|\la_{m+t}|\geq |\la_{m+t+l}|$ for all $l\in\{1,...,n-m-t\}$ and $\mu(\La)>\nu(4)$. From here $\mu(\La')\geq\mu(\La),$ so condition \textit{2} in the theorem is satisfied. Finally, to prove that condition \textit{3} in the theorem holds, we consider 
\begin{equation*}
N(b_{t+1})=N(\lambda_{m+t})^2N(a_{m+t})N\left(1+ (\la_{m+t}a_{m+t}\overline{\la_{m+t}})^{-1}\sum_{u=m+t+1}^n \la_u a_u\overline{\la_u}\right),
\end{equation*}
where $|(\la_{m+t}a_{m+t}\overline{\la_{m+t}})^{-1}|=|a_{m+t}|^{-1}|\la_{m+t}|^{-2}$. Since $|a_{m+t+l}|<|4a_{m+t}|$ and $|\la_{m+t}|\geq |\la_{m+t+l}|$ for all $l\in\{1,...,n-m-t\}$, we obtain that $$N\left(1+ (\la_{m+t}a_{m+t}\overline{\la_{m+t}})^{-1}\sum_{u=m+t+1}^n \la_ua_u\overline{\la_u}\right)$$ is a square due to Lemma \ref{LTCL}. We conclude that $N(b_{t+1})\in N(a_{m+t})k^{*2}$ and the proof of the condition \textit{3} is completed.
\end{dem1}

The following result, together with Lemma \ref{gen0}, give us an easy method to construct simple rotations of type (B) for binary lattices, from a given quaternion $r\in\OO_D$ as in Definition \ref{defAB}.

\begin{lema}
Let $r\in\OO_D$ be a non-zero quaternion and let $a_1,a_2\in\OO_D$ be non-zero pure quaternions. There exists $\la\in\OO_D$ different from zero such that $a_1=ra_1\barr+\la a_2\barl$ if and only if $NzNa_2\in k^{*2}$ and $NzNa_2^{-1}\in\OO_k$, where $z=a_1-ra_1\barr.$
\label{ctsu}
\end{lema} 

\begin{dem} The equation $a_1=ra_1\barr+\la a_2\barl$ has a solution $\la\in D^*$ if and only if the skew-hermitian form whose Gram matrix is $\left(\begin{array}{cc} z & 0 \\ 0 & -a_2 \end{array}\right)$ is isotropic. In fact, for $x,y\in D^*$, $xz\barx-ya_2\bary=0$ if and only if $z=\la a_2\barl$, where $\la=x^{-1}y$. Now, such a skew-hermitian form is isotropic if and only if it has discriminant 1 \cite[Ch. 10, \S3, Theorem 3.6]{Sch}. The discriminant is the reduced norm of its Gram matrix. We conclude that there exists $\la\in D^*$ such that $a_1=ra_1\barr+\la a_2\barl$ if and only if $NzNa_2\in k^{*2}.$ Finally, we have $Nz=Na_2N\la^2$ if $z=\la a_2\barl$. Hence, $\la\in\OO_D$ if and only if $NzNa_2^{-1}\in\OO_k.$
\end{dem}


We are ready to prove the Theorem \ref{eq}.

\begin{dem2} If $H(\Lambda)=k^*$, then there exists a simple rotation $(s;\s)\in\B(\La)$ such that $\theta[(s;\s)]\notin N\big(k(a_1)^*\big)/k^{*2}$ in virtue of Remark \ref{obscorkest}, so that 1) implies 2). To prove that 2) implies 3), let $(s;\s)$ be a simple rotation such that $\theta[(s;\s)]\notin N\big(k(a_1)^*\big)/k^{*2}$. As an isometry, $(s;\s)$ satisfies $a_1=h(s_1,s_1)=h((s;\s)(s_1),(s;\s)(s_1))=ra_1\barr+\lambda a_2\barl$ for some $r,\la\in\OO_D.$ Such an $r\in\OO_D$ satisfies $\s=a_1(1-\barr)$ (Lemma \ref{gen0}). Hence, $\theta[(s;\s)]\notin N\big(k(a_1)^*\big)/k^{*2}$ if and only if $N(1-r)\notin N\big(k(a_1)^*\big)$, in virtue of the equation $\theta[(s;\s)]=N(\s)k^{*2}$. If we rewrite the condition $N(1-r)\notin N\big(k(a_1)^*\big)$ in terms of the Hilbert symbol, we have $\left(\frac{N(1-r)\,,\,-Na_1}{\pe}\right)=-1$ since $a_1^2=-Na_1$. On the other hand, Lemma \ref{ctsu} tell us that $NzNa_2\in k^{*2}$ and $NzNa_2^{-1}\in\OO_k$, where $z=a_1-ra_1\barr.$ The result follows since $Na_2\in N(a_1)k^{*2}$ and $\mu=\nu(a_2)-\nu(a_1)=\nu(\pi^t)$. Finally, if $r\in\OO_D$ satisfies $NzNa_1\in k^{*2}$ and $NzN(\pi^ta_1)^{-1}\in \OO_k$, where $z=a_1-ra_1\barr$ and $\mu=\nu(\pi^t),$ then Lemma \ref{ctsu} and Lemma \ref{gen0} imply the existence of a simple rotation $(s;\s)\in\U_k^+(\La)$ of type (B) such that $\s=a_1(1-\barr).$ The condition on the Hilbert symbol implies that $\theta[(s;\s)]=N(\s)k^{*2}\notin N\big(k(a_1)^*\big)/k^{*2}$. Therefore $H(\Lambda)=k^*$ in virtue of Corollary \ref{kest}. This concludes the proof.
\end{dem2}

\begin{cor}
Let $\La$ be as in Theorem \ref{eq}. Let $t$ be such that $\mu=\nu(\pi^t)$. If $H(\La)=k^*$, then $H(\La')=k^*$ for every lattice $\La'=\langle b_1\rangle\bot\langle b_2\rangle$ with $N(b_1),N(b_2)\in N(a_1)k^{*2}$ and $\mu(\La')=\nu(\pi^s)$, for $s<t.$
\label{corteo}
\end{cor}

\begin{remark}
Due to Lemma \ref{v16}, in the condition \textit{2} of Theorem \ref{eq}, is enough to consider simple rotations $(s;\s)\in\B(\La)$ with $|\la|>|4|$, where $s=(1-r)s_1-\la s_2.$ Remember that $|1-r|\geq |2|$ for $(s;\s)\in\B(\La)$. Furthermore, we can assume that $|a_1|\geq |i|$ by rescaling. We use these facts in the proof of Theorem \ref{ssigen}.
\label{rla4}
\end{remark}

\begin{dem3}
If $r\in\OO_D$ satisfies the $k$-star conditions, then there exists $(s;\s)\in\B(\La)$, i.e., $s=(1-r_0)s_1-\la s_2$, with $|1-r_0|\geq 2$ and $\s=a_1(1-\overline{r_0})$, such that $N(1-r_0)\notin N\big(k(a_1)^*\big)$. Such an $r_0\in\OO_D$ also satisfies the $k$-star conditions. Let $\al\in\OO_D$ be a representative of the class of $r_0$ modulo $\pi^u$ as in the statement. Then, $r_0=\al+\pi^u\be$, with $\be\in\OO_D$ and $\al\in\SSS\oplus\SSS\om\oplus\SSS i\oplus\SSS i\om\subset\OO_D$. As $1-r_0=1-\al-\pi^u\be$ we have $N(1-r_0)=N(1-\al)N(1-(1-\al)^{-1}\pi^u\be)$. Now, $|1-r_0|\geq|2|$ implies $|1-\al|\geq|2|$. Therefore, $N(1-(1-\al)^{-1}\pi^u\be)$ is a square in virtue of Lemma \ref{LTCL}. Hence, $$\left(\frac{N(1-r_0)\,,\,-Na_1}{\pe}\right)=\left(\frac{N(1-\al)\,,\,-Na_1}{\pe}\right).$$ 
On the other hand, if $z=a_1-r_0a_1\overline{r_0}$ and $z'=a_1-\al a_1\overline{\al}$, then $z=z'-\pi^u\gamma$, with $\gamma=\al a_1\bar{\be}+\be a_1\bar{\al}+\pi^u\be a_1\bar{\be}\in\OO_D.$ Note that $a_1^{-1}\gamma\in\OO_D.$ We have $|z'|=|z|,$ since $|z|=|\pi^t\la a_1\bar{\la}|>|16\pi^ta_1|=|\pi^{4e+t}a_1|>|\pi^u\gamma|$, where we are assuming $|\la|>|4|$ (see Remark \ref{rla4}). Furthermore, we have that $Nz=Nz'N(1-z'^{-1}\pi^u\gamma)$ with $|z'^{-1}\pi^u\gamma|<|\pi^{-(4e+t)}a_1^{-1}\pi^{t+6e}\gamma|=|\pi^{2e}a_1^{-1}\gamma|\leq|4|.$ Hence, $NzNa_1 \textup{ is a square }$ if and only if $Nz'Na_1$ is a square. Finally, from $|z'|=|z|$, i.e., $|Nz|_k=|Nz'|_k$ we obtain $|Nz/\pi^{2t}N(a_1)|_k\leq 1$ if and only if $|Nz'/\pi^{2t}N(a_1)|_k\leq 1.$
\end{dem3}                                                                                                                                                                                                                                                                                                                                                                                                                                                                                                                                                                                                                                                                                                                                                                                                                                                                                                                                                                                                                                                                                                                                                                                                                                                                                                                                                                                                                                                                                                                                                                                                                                                                                                                                                                                                                                                                                                                                                                                                                                                                                             
                                                                                                                                                                                                                                                                                                                                                                                                                                                                                                                                                                                                                                                                                                                                                                                                                                                                                                                                                                                                                                                                                                                                                                                                                                                                                                                                                                                                                                                                                                                                                                                                                                                                                                                                                                                                                                                                                                                                                                                                                                                                                                                                                   
\begin{remark}
The number $u$ in Theorem \ref{ssigen} depends on $|\la|$. For example, if $\la$ satisfies $|\la|=1$, then we would have $|z'|=|\pi^t\la a_1\bar{\la}|=|\pi^ta_1|$ and so $|z'^{-1}\pi^u\gamma|=|\pi^{-t}a_1^{-1}\pi^u\gamma|\leq |\pi^{u-t}|<|4|$ if $u=t+3e.$ This holds in some cases when $k=\Q_2.$
\label{u}
\end{remark}                                                                                                                                                                                                                                                                                                                                                                                                                                                                                                                                                                                                                                                                                                                                                                                                                                                                                                                                                                                                                                                                                                                                                                                                                                                                                                                                                                                                                                                                                                                                                                                                                                                                                                                                                                                                                                                                                                                                                                                                                                                                                                                                                   
                                                                                                                                                                                                                                                                                                                                                                                                                                                                                                                                                                                                                                                                                                                                                                                                                                                                                                                                                                                                                                                                                                                                                                                                                                                                                                                                                                                                                                                                                                                                                                                                                                                                                                                                                                                                                                                                                                                                                                                                                                                                                                                                                                                                                                                                                                                                                                                                                                                                                                                                                                                                                                                                                                                                                                                                                                                                                                                                                                                                                                                                                                                                                                                                                                                                                                                                                                                                                                                                                                                                                                                                                                                                                                                                                                                                                                                                                                                                                                                                                                                                                                                                                                                                                                                                The following result let us choose particular lattices to compute $H(\Lambda)$ for arbitrary lattices. Note that for either of the remaining cases I or II described in the introduction, the extension $k(a_1)/k$ is ramified.

\begin{lema}
Let $\Lambda=\langle a_1\rangle \bot \langle a_2\rangle$ be a skew-hermitian lattice such that $N(a_2)\in N(a_1)k^{*2}$ and the extension $k(a_1)/k$ is ramified. Then, there exists a skew-hermitian lattice $L=\langle q\rangle\bot\langle \epsilon q\rangle$, where $q\in D^*$ and $\epsilon\in k^*,$ such that $H(\Lambda)=H(L)$. Moreover, we can assume that $q=q'$, for any quaternion $q'\in D^*$ with $N(q')\in N(a_1)k^{*2}.$
\label{ret}
\end{lema}                                                                                                                                                                                                                                                                                                                                                                                                                                                                                                                                                                                                                                                                                                                                                                                                                                                                                                                                                                                                                                                                                                                                                                                                                                                                                                                                                                                                                                                                                                                                                                                                                                                                                                                                                                                                                                                                                                                                                                                                                                                                                                                                                                                                                                                                                                                                                                                                                                

\begin{dem}
If $N(a_2)=N(a_1)b^2=N(ba_1)$ for some $b\in k^*$ and $k(a_1)/k$ is ramified, we can use Lemma 4.3 in \cite{A-C04} to conclude that $\langle a_1\rangle \bot \langle a_2\rangle\cong \langle a_1\rangle \bot \langle \xi ba_1\rangle$, some $\xi\in k^*.$ Therefore, it is enough to take $L=\langle q\rangle\bot\langle \epsilon q\rangle$, where $q=a_1$ and $\epsilon=\xi b$. Finally, if $q'\in D^*$ with $N(q')\in N(a_1)k^{*2}$, we have $N(q')=N(ca_1)$ for some $c\in k^*$. Applying Lemma 4.3 in \cite{A-C04} again, with $q=a_1$, we get $q'=\eta \varepsilon cq\bar{\eta}$, for some $\varepsilon\in k^*$ and $\eta\in D^*.$ Hence, $\Lambda'=\langle q'\rangle\bot\langle \epsilon q'\rangle$ is isometric to a rescaling of $L$. We conclude $H(\Lambda')=H(L)$ as stated.
\end{dem}

\section{Algorithm for $k=\Q_2$ and proof of Theorem \ref{completa}}\label{algq2}

By considering Theorems \ref{eq} and \ref{ssigen}, we are in conditions to construct an algorithm to compute $H(\La)$, for a binary lattice $\La$, as follows: If $k=\Q_2$, we have $\OO_k=\Z_2$, $\OO_D=\Z_2\oplus\Z_2\om\oplus\Z_2i\oplus\Z_2i\om$ and $\OO_k/\pi^u\OO_k\cong \Z/2^u\Z$.

\begin{enumerate}
\item Replace $\La$ by a suitable lattice $L=\langle q\rangle\bot\langle \epsilon q\rangle$ as in Lemma \ref{ret}. 
We have two possibilities for the square class of $Na_1:$
\begin{enumerate}
\item If $Na_1\in\al\Q_2^{*2}$, where $\al$ is a unit, we look for a pure quaternion $q\in\OO_D$ such that $Nq\equiv \al\;(8)$. In this case, $Nq=\al(1+8\al^{-1}\be)$, for some $\be\in\OO_D$. Hence $Nq\in N(a_1)\Q_2^{*2}$ and $|q|=1$.
\item If $Na_1\in\al\Q_2^{*2}$, where $\al$ is a prime, we look for a pure quaternion $q\in\OO_D$ such that $Nq\equiv \al\;(16)$. In this case, $Nq=\al(1+16\al^{-1}\be)$, for some $\be\in\OO_D$. Hence, $Nq\in N(a_1)\Q_2^{*2}$ and $|q|=|i|.$
\end{enumerate}
In both cases we obtain a pure quaternion $q\in\OO_D$ with $|q|\geq|i|.$

\item Fix a set of representatives $\SSS$ of the finite ring $\OO_k/\pi^u\OO_k:$ We can choose $\SSS=\{0,1,...,2^u-1\}$ because $\OO_k/\pi^u\OO_k\cong \Z/2^u\Z.$

\item For $r=a+b\om+ci+di\om\in\SSS\oplus\SSS\om\oplus\SSS i\oplus\SSS i\om\subset\OO_D$, check if the $k$-star conditions are satisfied. This verification can be done by using the Sage functions as $\al\mapsto (\al).ordp()$, which give us the $p$-adic valuation of $\al$ in $\Q_p$ \cite{S}. We know that $N(a,b,c,d,\pi)=a^2+ab-b^2-\pi(c^2+cd-d^2)$ is the norm of $r=a+b\om+ci+di\om\in\OO_D$. Then, if we write $z=q-rq\barr=z_0+z_1\om+z_2i+z_3i\om$, we get $Nz=N(z_0,z_1,z_2,z_3,\pi)$. Now, we write a little program using Sage as shown in the following example. Note that $Nz.subs()$ let us substitute values in the expression for $Nz$:



\begin{algorithm}
\caption{Case I with $\La=\langle j+ij\rangle\bot\langle 4(j+ij)\rangle$}
 R=Qp(2);
\begin{algorithmic}
\FOR{$a,b,c,d\in\SSS$}
   \IF{hilbert\_symbol$(N(a-1,b,c,d,2),-5,2)==-1$ and R$(Nz.subs(a=a,b=b,c=c,d=d)*5)$.is\_square() and R$(Nz.subs(a=a,b=b,c=c,d=d)/2^4).$ordp()$>=0$}
            \RETURN $H(\La)=\Q_2^*$
  \ENDIF
\ENDFOR
\RETURN $H(\La)=N\big(\Q_2(q)^*\big)$
\end{algorithmic}
\end{algorithm}

\item Conclude that $H(\La)=\Q_2^*$ if some $r$ in the last step satisfy the $k$-star conditions. Otherwise, conclude that $H(\La)=N\big(\Q_2(a_1)^*\big)$ in virtue of Theorem \ref{ssigen} and Proposition \ref{konok}.
\end{enumerate}                                                                                                                                                                                                                                                                                                                                                                                                                                                                                                                                                                                                                                                                                                                                                                                                                                                                                                                                                                                                                                                                                                                                                                                                                                                                                                                                                                                                                                                                                                                                                                                                                                                                                                                                                                                                                                                                                                                                                                                                                                                                                                                                                                                                                                                                                                                                                                                                                                                                                                                                                                                                                                                                                                                                                                                                                                                                                                                                                                                                                                                                                                                                                                                                                                                                                                                                                                                                                                                                                                                                                                                                                                                                                                                                                                                                                                                                                                                                                                                                                                                                                                                                                                                                                                                 
          
\begin{remark}
The algorithm depends only on $a_1$ and $\mu(\La)$. Then, for the lattice $L$ above, we are interested only in $q$ and the valuation of $\epsilon$, so we can multiply $\epsilon$ for any element of $\OO_k^*$.
\label{multk}
\end{remark}
                                                      
\begin{remark}
The algorithm can be extended to any unramified finite extension $k$ of $\Q_2$. The condition $|2a_1|\geq |a_2|$ in Theorems \ref{eq} and \ref{red} is essential. Hence, the algorithm does not work, for $\mu<\nu(2)$, if the extension $k/\Q_2$ ramifies, unless the algorithm returns the value $k^*$ for $\mu<\nu(2).$                                                                                                                                                                                                                                                                                                                                                                                                                                                                                                                                                                                                                                                                                                                                                                                                                                                                                                                                                                                                                                                                                                                                                                                                                                                                                                                                                                                                                                                                                                                                                                                                                                                                                                                                                                                                                                                                                                                                                                                                                                                                                                                                                 \end{remark}                                                                                                                                                                                                                                                                                                                                                                                                                                                                                                                                                                                                                                                                                                                                                                                                                                                                                                                                                                                                                                                                                                                                                                                                                                                                                                                                                                                                                                                                                                                                                                                                                                                                                                                                                                                                                                                                                


\subsection{Computations using Sage} 

To compute the spinor images in cases I and II when $k=\Q_2$, we use the algorithm above. The following results are obtained by computer search. When the algorithm actually find solutions, we actually list them. Otherwise it is just stated that no solutions were found.

\begin{lema}
For any $q\in\{j+ij,i+j\}$ and $t\in\{3,4\}$, there exist $r_1,r_2\in\OO_D$ such that:
\begin{enumerate}
\item $|1-r_1|=|2|,$ $NzNq\in\Q_2^{*2}$ and $NzN(2^tq)^{-1}\in\Z_2^*,$ where $z=q-r_1q\overline{r_1}$.
\item $|1-r_2|=|i|,$ $NzNq\in\Q_2^{*2}$ and $NzN(2^tq)^{-1}\in\Z_2^*,$ where $z=q-r_2q\overline{r_2}$.
\end{enumerate}

\label{v1}
\end{lema}

\begin{dem}
It is a direct computation to verify that the elements $r_1,r_2\in\OO_D$ in the table below satisfy the conditions \textit{1-2} in the lemma.

\begin{table}[h]
\begin{center}
    $$\begin{array}{|c|c|c|c|}
 \hline   
   q      & t  &r_1 & r_2\\ \hline
   j+ij    & 3  &-1-4i-4i\om&1-14\om-i-10i\om\\ \hline
   j+ij   & 4  &-1-8i-8i\om&1-6\om-13i-6i\om\\ \hline
   i+j    & 3  &-1-4i\om& 1-2\om-i\\ \hline
   i+j    & 4  &-1-8i\om& -1-6\om-3i\\ \hline
    \end{array}$$
\end{center}
\end{table}
\end{dem}

Using the Sage algorithm above we find elements $r\in\OO_D$ (see two following lemmas) which satisfy the $k$-star conditions for particular lattices. This help us to conclude in next section that $H(\La)=\Q_2^*$ in some cases. 

\begin{lema}
Let $\La=\langle a_1\rangle\bot\langle a_2\rangle$ be a skew-hermitian lattice satisfying the conditions in Theorem \ref{eq}. For $a_1\in\{j+ij,i+j\}$ and $t\in\{1,2\}$ there exists $r\in\OO_D$ satisfying the $k$-star conditions.
\label{kestu}
\end{lema}

\begin{dem}
It is a direct computation to verify that the elements $r\in\OO_D$ shown in the table below, satisfy the required conditions for $t=2$ and then, for $t=1$ in virtue of Corollary \ref{corteo}. Here, $i_{_2}=i$ satisfies $i_{_2}^2=2.$

\begin{center}
\begin{footnotesize}
\begin{tabular}{|c|c|c|c|c|}\hline
              
$a_1$ & $r$ & $N(1-r)$ & $z$ & $NzNa_1$\\ \hline
$j+ij$  & $1+2i_{_2}\om$ & $2\cdot 2^2$ & $4(-1+2\om-4i_{_2}-7i_{_2}\om)$ & $2^4\cdot 5^2$ \\ \hline 
$i+j$  & $1+2i_{_2}+2i_{_2}\om$ & $-2\cdot 2^2$ & $4(1-2\om+3i_{_2}+3i_{_2}\om)$ & $2^4(1+8\cdot 20)$ \\ \hline
\end{tabular}
\end{footnotesize}
\end{center}

\end{dem}

Note that $D=\left(\frac{\pi,\Delta}{k}\right)$ for any prime $\pi$ of $k$. Hence, for every prime $\pi$, there exists a pure quaternion $i_{_\pi}\in ik(j)$ satisfying $i_{_\pi}^2=\pi$ and $i_{_\pi}j=-i_{_\pi}j.$                                                                                                                                                                                                                                                                                                                                                                                                                                                                                                                                                                                                                                                                                                                                                                                                                                                                                                                                                                                                                                                                                                                                                                                                                                                                                                                                                                                                                                                                                                                                                                                                                                                                                                                                                                                                                                                                                                                                                                                                                                                                                                                                                                                                                                                                                                                                                                                                                                                                                                                                 

\begin{lema}
Let $\La=\langle a_1\rangle\bot\langle a_2\rangle$ be a skew-hermitian lattice satisfying the hypothesis in Theorem \ref{eq}. For every $a_1\in\{i_{_\pm 2},i_{_\pm 10}\}$ as above and $t\in\{1,2,3,4\}$, there exists $r\in\OO_D$ satisfying the $k$-star conditions.
\label{kestpi}
\end{lema}

\begin{dem}
As in Lemma \ref{kestu}, it is easy to see that the elements $r\in\OO_D$ shown in the table below satisfy the required conditions for $t=4$ and then, for $t<4$ in virtue of Corollary \ref{corteo}.

\begin{center}
\begin{footnotesize}
\begin{tabular}{|c|c|c|c|c|}\hline
              
$\pi$ & $r$ & $N(1-r)$ & $z$ & $NzNa_1$\\ \hline
$\pm 2$ & $15+8\om$ & $5\cdot 2^2(1+8\cdot 5^{-1}\cdot 7)$ & $-592i_{\pi}+304i_{_\pi}\om$ & $2^{10}(1+8\cdot 38)$\\ \hline
$\pm 10$ & $15+8\om$ & $5\cdot 2^2(1+8\cdot 5^{-1}\cdot 7)$ & $-592i_{\pi}+304i_{_\pi}\om$ & $2^{10}\cdot 5^2(1+8\cdot 38)$\\ \hline       
\end{tabular}
\end{footnotesize}
\end{center}                                                                                                                                                                                                                                                                                                                                                                                                                                                                                                                                                                                                                                                                                                                                                                                                                                                                                                                                                                                                                                                                                                                                                                                                                                                                                                                                                                                                                                                                                                                                                                                                                                                                                                                                                                                                                                                                                                                                                                                                                                                                                                                                                                                                                                                                                                                                                                                                                                                                                                                                                                                                                                                                                                                                                                                                                                                                                                                                                                                                                                                                                                                                                                                                                                                                                                                                                                                                                                                                                                                                                                                                                                                                                                                                                                                                                                                                                                                                                                                                                                                                                                                                                                                                                                                   

\end{dem}

The following result help us to prove, in next section, that $H(\La)\neq \Q_2^*$ in some cases. It is proved by an explicit search using Sage as above.

\begin{lema}
There is no $r=a+b\om+ci+di\om\in \Z\oplus\Z\om\oplus\Z i\oplus\Z i\om\subset\OO_D$, with $0\leq a,b,d,c<2^{t+3}$ satisfying the $k$-star conditions for $t\in\{3,4\}$ and $a_1\in\{j+ij,j+i\}.$
\label{noexistev8}
\end{lema}
                                                                                                                                                                                                                                                                                                                                                                                                                                                                                                                                                                                                                                                                                                                                                                                                                                                                                                                                                                                                                                                                                                                                                                                                                                                                                                                                                                                                                                                                                                                                                                                                                                                                                                                                                                                                                                                                                                                                                                                                                                                                                                                                                                                                                                                                                                                                                                                                                                                                                                                                                                                                                                                                                                                                                                                                                                                                                                                                                                                                             

\subsection{Proof of Theorem 1 in Cases I and II}

                                                                                                                                                                                                                                                                                                                                                                                                                                                                                                                                                                                                                                                                                                                                                                                                                                                                                                                                                                                                                                                                                                                                                                                                                                                                                                                                                                                                                                                                                                                                                                                                                                                                                                                                                                                                                                                                                                                                                                                                                                                                                                                                                                                                                                                                                                                                                                                                                                                                                                                                                                                                                                                                  \textbf{Proof in Case I.} Here we have $\Lambda=\langle a_1\rangle\bot ...\bot\langle a_n\rangle$, where $N(a_m)\in-u\Q_2^{*2},$ for each $m=1,...,n$ and $u\in\Z_2^*$ is a unit of non-minimal quadratic defect independent of $m$. As $\Z_2^*/\Z_2^{*2}=\{\overline{\pm 1},\overline{\pm 5}\}$ and a pure quaternion cannot have norm $-1$, we have two options for $u$: $u=-5$ or $u=-1$.

In virtue of Lemma \ref{ret}, we consider for binary lattices, $\Lambda=\langle q\rangle\bot\langle \epsilon q\rangle$, where we can choose any pure quaternion $q\in\OO_D^0$ satisfying $N(q)\in -u\Q_2^{*2}$, and $\epsilon=\varepsilon\pi^t$, with  $\varepsilon\in\Z_2^*$ and $t\in\{1,2,3,4\}$. Here, $q=q_u$ runs over a system of representatives with $N(q)\in-u\Q_2^{*2}$, for $u$ running over the set $\{-5,-1\}$ of units of non-minimal quadratic defect. Moreover, we can assume $i=i_{_2}$, so $i^2=2$, and therefore 
we work with $\Lambda=\langle q_u\rangle\bot\langle 2^tq_u\rangle$, where $t\in\{1,2,3,4\},$ $q_{_{-5}}=j+ij$ and $q_{_{-1}}=i+j$. 

\begin{prop}
Let $\Lambda=\langle a_1\rangle\bot ...\bot\langle a_n\rangle$ be a skew-hermitian lattice such that $N(a_1)$,..., $N(a_n)\in -u\Q_2^{*2}$ and $0<\mu(\La)<\nu(8)$. Then $H(\Lambda)=\Q_2^*.$                                                                                                                                                                                                                                                                                                                                                                                                                                                                                                                                                                                                                                                                                                                                                                                                                                                                                                                                                                                                                                                                                                                                                                                                                                                                                                                                                                                                                                                                                                                                                                                                                                                                                                                                                                                                                                                                                                                                                                                                                                                                                                                                   
\label{propcaso1k}
\end{prop}                                                                                                                                                                                                                                                                                                                                                                                                                                                                                                                                                                                                                                                                                                                                                                                                                                                                                                                                                                                                                                                                                                                                                                                                                                                                                                                                                                                                                                                                                                                                                                                                                                                                                                                                                                                                                                                                                                                                                                                                                                                                                                                    

\begin{dem}
It is enough to consider the case $n=2$. So by the discussion above, we can assume that $\Lambda=\langle q_u\rangle\bot\langle 2^tq_u\rangle,$ with $q_{_{-5}}=j+ij$ and $q_{_{-1}}=i+j$, where $t\in\{1,2\}.$ In virtue of Corollary \ref{corteo} it suffices to prove the result for $t=2.$ Lemma \ref{kestu} tell us that there exists $r\in\OO_D$ satisfying the $k$-star conditions. This is equivalent to $H(\La)=\Q_2^*.$
\end{dem}

To handle the cases where $\nu(8)\leq \mu\leq\nu(16)$ we use the following result, which is used to improve the set of generators $\B(\La).$ The proof is a routine calculation.

\begin{lema}                                                                                                                                                                                                                                                                                                                                                                                                                                                                                                                                                                                                                                                                                                                                                                                                                                                                                                                                                                                                                                                                                                                                                                                                                                                                                                                                                                                                                                                                                                                                                                                                                                                                                                                                                                                                                                                                                                                                                                                                                                                                                                                       
If $r\in \OO_D$ satisfies any of the equations 
\begin{eqnarray*}
j+ij & = & r(j+ij)\barr+\ep 2^t\la (j+ij)\barl, \\
i+j & = & r(i+j)\barr+\ep 2^t\la (i+j)\barl,
\label{rjkrc}
\end{eqnarray*}
where $\ep\in\Z_2^*,\;\la\in\OO_D,\;t\geq 2$, then $1-r\in i\OO_D.$                                                                                                                                                                                                                                                                                                                                                                                                                                                                                                                                                                                                                                                                                                                                                                                                                                                                                                                                                                                                                                                                                                                                                                                                                                                                                                                                                                                                                                                                                                                                                                                                                                                                                                                                                                                                                                                                                                                                                                                                                                                                                                                   
\label{lrjkrc}
\end{lema}

\begin{lema}
Let $\Lambda=\OO_Ds_1\bot\OO_Ds_2=\langle a_1\rangle\bot\langle a_2\rangle$ be a skew-hermitian lattice such that $N(a_1),N(a_2)\in-u\Q_2^{*2}$ 
and $\nu(8)\leq \mu(\La)\leq\nu(16)$. There exists a lattice $L$ of rank 2 such that $H(L)=H(\La)$ and the subsets $\B_1(L),\B_2(L)\subset\B(L)$, defined by $\B_1(L)=\{(s;\s)\in\B(L)\;:\;|1-r|=|i|\}$ and $\B_2(L)=\{(s;\s)\in\B(L)\;:\;|\la|=1\}$, where $r$ as in Lemma \ref{gen0} and $\la$ as in Lemma \ref{ctsu}, satisfy $\mathcal{A}(L)\cup \B_l(L)$ generates $\U_{_{\Q_2}}^+(L),$ for $l=1,2.$
\label{genq2}
\end{lema}

\begin{dem} Remember that $u\in\{-5,-1\}.$ By Lemma \ref{ret} the lattice $L=\OO_Ds_1\bot\OO_Ds_2=\langle q_u\rangle\bot\langle \ep 2^tq_u\rangle,$ where $\ep\in\Z_2^*$ and $t\in\{3,4\}$ satisfies $H(L)=H(\La)$, where $q_{_{-5}}=j+ij$ and $q_{_{-1}}=i+j$. Let $\phi\in\B(L)$ be such that $\phi(s_1)=rs_1+\la s_2$. We have $|1-r|\in\{|i|,|2|\}$ in virtue of Lemma \ref{lrjkrc}. Hence, to prove that $\B_1(L)$ satisfies the required property, it suffices to prove\footnote{In this case, there exists a second element $(s';\s')\in\B_1(L)$ defined by $s'=s_1-(s;\s)\phi(s_1),\;\s'=q(1-\overline{r''})$ such that $(s';\s')(s;\s)\phi(s_1)=s_1$.} that, if $\phi$ satisfies $|1-r|=|2|$, then there exists $(s;\s)\in\B(L)$ such that $|1-r''|=|i|$ and $|1-r'|=|i|$, where $r',r''\in\OO_D$ are defined by $(s;\s)(s_1)=r's_1+\la's_2$ and $(s;\s)\phi(s_1)=r''s_1+\la''s_2$, for some $\la',\la''\in\OO_D$. In fact, computing $(s;\s)\phi(s_1)=(s;\s)(rs_1+\la s_2)$ we have 

\begin{eqnarray}
1-r'' & = & 1-r+[rq(1-\bar{r'})+\la \ep 2^tq\bar{\la'}](1-\bar{r'})^{-1}q^{-1}(1-r')\label{1-r''},\\  
\la'' & = &\la+[rq(1-\bar{r'})+\la \ep 2^tq\bar{\la'}](1-\bar{r'})^{-1}q^{-1}\la'.\label{la''}
\end{eqnarray}

Lemma \ref{v1} implies the existence of an element $r'\in\OO_D$ 
 such that $$|1-r'|=|i|,\hspace{1cm} NzN(\ep 2^tq)\in\Q_2^{*2} \hspace{0.5cm}\text{  and  }\hspace{0.5cm} NzN(\ep 2^tq)^{-1}\in\Z_2,$$ where $z=q-r'q\overline{r'}$. Hence, by Lemma \ref{ctsu}, there exists $\la'\in\OO_D$, such that $q=r'q\overline{r'}+\ep 2^t\la'q\overline{\la'}$. Then $(s;\s)$ with $s=(1-r')s_1-\la's_2,\;\s=q(1-\bar{r'})$ defines a simple rotation of type (B) in virtue of Lemma \ref{gen0}. Note that $(s;\s)\in\B_1(L)$. On the other hand, as 

\begin{eqnarray*}
|[rq(1-\bar{r'})+\la \ep 2^tq\bar{\la'}](1-\bar{r'})^{-1}q^{-1}(1-r')|&=&|rq(1-\bar{r'})+\la \ep 2^tq\bar{\la'}|\\
&=&|1-\bar{r'}|=|i|
\end{eqnarray*}
 and $|1-r|=|2|$, it follows that $|1-r''|=|i|.$ In particular, $\mathcal{A}(L)\cup \B_1(L)$ generates $\U_{_{\Q_2}}^+(L).$
 

Now, to prove that $\mathcal{A}(L)\cup\B_2(L)$ generates $\U^+(L)$, it suffices to prove\footnote{By a similar argument as we did for $\B_1(L).$ 
} that, if $\phi\in\B(L)$ satisfies $|\la|<1$, there exists $(s;\s)\in\B(L)$ such that $|\la'|=1$ and $|\la''|=1$, where $\la,\la',\la''$ are defined by $\phi,(s;\s)$ and $(s;\s)\phi$ respectively, as before. From the equation (\ref{la''}) we see that $|\la''|=1$ if $|\la|<1$ and $|\la'|=1.$ 
By Lemma \ref{v1}, there exists $r'\in\OO_D$  
such that $$|1-r'|=|i| \text{ or } |2|, \hspace{0.5cm} NzN(\ep 2^tq)\in\Q_2^{*2} \hspace{0.5cm}\text{  and  }\hspace{0.5cm} NzN(\ep2^tq)^{-1}\in\Z_2^*,$$ where $z=q-r'q\overline{r'}$. Hence, by Lemma \ref{ctsu}, there exists $\la'\in\OO_D$ such that $q=r'q\overline{r'}+\ep 2^t\la'q\overline{\la'}$. Then $(s;\s)$, with $s=(1-r')s_1-\la's_2$ and $\s=q(1-\bar{r'})$, defines a simple rotation of type (B) in virtue of Lemma \ref{gen0}, and this rotation satisfies $|\la'|=1$ since $NzN(\ep2^tq)^{-1}\in\Z_2^*.$ Now, we take $|1-r'|=|i|, \text{ if } |1-r|=|2|, \text{ and } |1-r'|=|2|,\text{ if } |1-r|=|i|,$ so that $|1-r'|,|1-r''|\geq |2|$ in virtue of equation (\ref{1-r''}). The result follows.
 
 
\end{dem}

\begin{remark}
Notice that, for a lattice $\La$ as in the previous lemma, we can replace $\B(\La)$ by $\B_l(\La), $ for $l=1,2$, in Theorem \ref{eq}. Hence, for any $(s;\s)$ in $\B_2(L)$ we have that, if $z=q-rq\barr$, then $|z|=|2^t|$ since $|z|=|2^t\la q\bar{\la}|.$ This fact, help us to improve the number $u$ in Theorem \ref{ssigen} in virtue of Remark \ref{u}. 
\label{BporB1B2}
\end{remark}


\begin{prop}
There exists $r\in \OO_D$ satisfying the $k$-star conditions for $t\in\{3,4\}$ and $Na_1\in -u\Q_2^{*2}$, with $u$ a unit of non-minimal quadratic defect if and only if there exists $\al=a+b\om+ci+di\om\in \Z\oplus\Z\om\oplus\Z i\oplus\Z i\om\subset\OO_D$, with $0\leq a,b,d,c<2^{t+3}$, satisfying them.
\label{ssi}
\end{prop}

We have a direct consequence of Proposition \ref{genq2} and Lemmas \ref{ret}, \ref{noexistev8}.

\begin{cor}
Let $\La=\langle a_1\rangle\bot\langle a_2\rangle$ be a skew-hermitian lattice such that $N(a_1),N(a_2)\in -u\Q_2^{*2}$, where $u$ is a unit of non-minimal quadratic defect and $\mu=\nu(a_2)-\nu(a_1)$ satisfies $\nu(8)\leq\mu\leq\nu(16)$. Then $H(\La)=N\big(\Q_2(a_1)^*\big).$                                                                                                                                                                                                                                                                                                                                                                                                                                                                                                                                                                                                                                                                                                                                                                                                                                                                                                                                                                                                                                                                                                                                                                                                                                                                                                                                                                                                                                                                                                                                                                                                                                                                                                                                                                                                                                                                                                                                                                                                                                                                                                                                   
\label{nq2v8}
\end{cor}

We need the following result to handle lattices $\La$ with $\mu(\La)=\nu(8).$

\begin{lema}
If $|\eta|=|i|$ and $a_1$ is a pure unit, then $T(2(\eta a_1\bar{\eta})^{-1}a_1)\in\pi\OO_k.$
\label{traza}
\end{lema}

\begin{dem}
Set $\rho=i^{-1}\eta\in\OO_D^*,$ so that $\eta=i\rho.$ Note that $a_1i\equiv i\bar{a_1}\text{ (mod }\pi),$ where $\rho$ and $\bar{a_1}$ commute modulo $i$. We conclude that $$\eta a_1\bar{\eta}\equiv -N(\rho)\pi \bar{a_1}\text{ (mod }\pi i).$$ In other words $\frac{1}{\pi}\eta a_1\bar{\eta}=-N(\rho)\bar{a_1}+\varepsilon,$ where $\varepsilon\in i\OO_D,$ whence $\pi(\eta a_1\bar{\eta})^{-1}=-(N(\rho)\bar{a_1})^{-1}+\delta=\frac{-a_1}{N(\rho a_1)}+\delta$, for some $\delta\in i\OO_D.$ We conclude that $$T\big(2(\eta a_1\bar{\eta})^{-1}a_1\big)\equiv \frac{-4a_1^2}{\pi N(\rho a_1)}+\frac{2}{\pi}T(\delta a_1)\text{ (mod }\pi)$$ and the result follows since $\delta\in i\OO_D$ implies $T(\delta)\in\pi\OO_k.$
\end{dem}

\begin{prop}
Let $\La=\langle a_1\rangle\bot ...\bot\langle a_n\rangle$ be a skew-hermitian lattice such that $N(a_1),...,$ $ N(a_n)\in -u\Q_2^{*2}$, where $u$ is a unit of non-minimal quadratic defect. If $\mu=\mu(\La)$ satisfies $\nu(8)\leq\mu\leq\nu(16)$, then $H(\La)=N\big(\Q_2(a_1)^*\big).$                                                                                                                                                                                                                                                                                                                                                                                                                                                                                                                                                                                                                                                                                                                                                                                                                                                                                                                                                                                                                                                                                                                                                                                                                                                                                                                                                                                                                                                                                                                                                                                                                                                                                                                                                                                                                                                                                                                                                                                                                                                                                                                                   
\label{propcaso1}
\end{prop}                                                                                                                                                                                                                                                                                                                                                                                                                                                                                                                                                                                                                                                                                                                                                                                                                                                                                                                                                                                                                                                                                                                                                                                                                                                                                                                                                                                                                                                                                                                                                                                                                                                                                                                                                                                                                                                                                                                                                                                                                                                                                                                    

\begin{dem} We consider separately the cases $\mu=\nu(16),\mu=\nu(8).$\\

\underline{\textbf{Case 1:} $\mu=\nu(16).$}\\

In virtue of Lemma \ref{v16} it suffices to consider rotations $(s;\s)\in\B(\La)$ such that $|1-r|=|2|$ and $|\la_2|=1$. In this case, Theorem \ref{red} tell us we can set $n=2$ in the statement of the proposition. For $n=2$, because of Lemma \ref{genq2}, we can replace $\La$ by a lattice $L$ such that $H(L)=H(\La)$ and a set of generators of $\U_{_{\Q_2}}^+(L)$ is $\mathcal{A}(L)\cup\B_1(L)$. It follows that $H(\Lambda)=N\big(\Q_2(a_1)^*\big)$ since rotations in $\B_1(L)$ have spinor norm belonging to $N\big(\Q_2(a_1)^*\big)$ in virtue of Lemma \ref{v16}. \\

\underline{\textbf{Case 2:} $\mu=\nu(8).$} \\

In virtue of Lemma \ref{v16}, any rotation $(s;\s)\in\B(\La)$ satisfy $\theta[(s;\s)]\in N\big(\Q_2(a_1)^*\big)/\Q_2^{*2}$ unless one of the following conditions is satisfied: 
\begin{enumerate}
\item $|1-r|=|i|,\;|\la_2|=1$,
\item $|1-r|=|2|,\;|\la_2|\in\{1,|i|\}.$
\end{enumerate}
As in the Case 1 we can reduce the cases for which $|\la_2|=1$ to consider rank 2 lattices and the case $|1-r|=|2|,\;|\la_2|=|i|$ to study rank 3 lattices with\footnote{If $|\la_3|<1$, then we can reduce the study to rank 2 lattices (see Theorem \ref{red}).} $|\la_3|=1$. For rank 2 lattices, Corollary \ref{nq2v8} tell us that $H(\Lambda)=N\big(\Q_2(a_1)^*\big).$ We prove that, for rank 3 lattices $\La$ such that $(s;\s)\in\B(\La)$ satisfies $|1-r|=|2|,\;|\la_2|=|i|,\;|\la_3|=1$ we also have $\theta[(s;\s)]\in N\big(\Q_2(a_1)^*\big)/\Q_2^{*2}.$ In fact, in virtue of Lemma \ref{ret} we can assume that $\La=\langle a_1\rangle\bot\langle 8\ep_2 a_1\rangle\bot\langle 64\ep_3 a_1\rangle$, with $\ep_2,\ep_3\in\Z_2^*$. Hence, Lemma \ref{gen0} tell us that $r,\la_2,\la_3\in\OO_D$, with $|1-r|\geq |2|$, define an element $\phi\in\B(\La)$ if and only if they satisfy the relation $$z=8\la_2 \ep_2a_1\overline{\la_2}+64\la_3\ep_3a_1\overline{\la_3},$$ where $z=a_1-ra_1\barr$. We can rewrite this equation as follows: $$z=8\la_3 (\ep_2\eta a_1\overline{\eta}+8\ep_3a_1)\overline{\la_3},$$ where $\eta=\la_3^{-1}\la_2.$ Remember that, in this case, $|\la_2|=|i|$ and $|\la_3|=1$. Hence, by Lemma \ref{ctsu}, the existence of $r,\la_1,\la_2$ satisfying the equation above is equivalent to the existence of $r,\eta\in\OO_D$, with $|\eta|=|i|$ such that $NzN(\ep_2\eta a_1\overline{\eta}+8\ep_3 a_1)\in\Q_2^{*2}$ and $ NzN\big(8(\ep_2\eta a_1\overline{\eta}+8\ep_3 a_1)\big)^{-1}\in\Z_2$. We know that $|\ep_2\eta a_1\overline{\eta}+8\ep_3 a_1|=|2|$, so $NzN\big(8(\ep_2\eta a_1\overline{\eta}+8\ep_3 a_1)\big)^{-1}\in\Z_2$ if and only if $\frac{Nz}{2^8}\in\Z_2.$ On the other hand, $N(\ep_2\eta a_1\overline{\eta}+8\ep_3 a_1)=N(\ep_2\eta a_1\overline{\eta})N(1+8\ep(\eta a_1\overline{\eta})^{-1}a_1)$, where $\ep=\ep_2^{-1}\ep_3\in\Z_2^*$. Here, as $|8\ep(\eta a_1\overline{\eta})^{-1}a_1|=|4|$, we can write $8\ep(\eta a_1\overline{\eta})^{-1}a_1=4\ep\xi$ with $\xi=2(\eta a_1\overline{\eta})^{-1}a_1\in\OO_D^*$ and we have that $$N(1+4\ep\xi)=1+4\ep T(\xi)+16\ep^2N(\xi).$$ As $T(\xi)=2x+y$, with $x,y\in\Z_2$ it follows that $N(1+4\ep\xi)\in\{\Q_2^{*2},\;5\Q_2^{*2}\}.$ Hence, $N(\ep_2\eta a_1\overline{\eta}+8\ep_3 a_1)\in\{N(a_1)\Q_2^{*2},\;5N(a_1)\Q_2^{*2}\}$, where $N(\ep_2\eta a_1\overline{\eta}+8\ep_3 a_1)\in 5Na_1\Q_2^{*2}$ if and only if $T(\xi)\equiv 1$ (mod 2). The last condition is not satisfied in virtue of Lemma \ref{traza}. Therefore, we are reduced to the following result, that is an analogue of Theorem \ref{eq}: \textit{$H(\La)=\Q_2^*$ if and only if there exists $r\in\OO_D$, with $|1-r|=|2|$ such that:
\begin{enumerate}
\item $\big(\frac{N(1-r),-Na_1}{\pe}\big)=-1$,
\item $NzNa_1\in\Q_2^{*2}$,
\item $|z|=|16|$.
\end{enumerate}}
Corollary \ref{nq2v8} implies that there is no $r\in\OO_D$ satisfying the conditions above. Hence, we conclude that $H(\La)=N\big(\Q_2(a_1)^*\big).$

\end{dem}                                                                                                                                                                                                                                                                                                                                                                                                                                                                                                                                                                                                                                                                                                                                                                                                                                                                                                                                                                                                                                                                                                                                                                                                                                                                                                                                                                                                                                                                                                                                                                                                                                                                                                                                                                                                                                                                                                                                                                                                                                                                                                                                                                                                                                                                                                                                                                                                                                                                                                                                                                                                                                          
                                                                                                                                                                                                                                                                                                                                                                                                                                                                                                                                                                                                                                                                                                                                                                                                                                                                                                                                                                                                                                                                                                                                                                                                                                                                                                                                                                                                                                                                                                                                                                                                                                                                                                                                                                                                                                                                                                                                                                                                                                                                                                                                      
\textbf{Proof in Case II.} The following result is valid over an arbitrary dyadic local field $k$ and it is an improvement of Proposition 6.9 in \cite{A-C04}.                                                                                                                                                                                                                                                                                                                                                                                                                                                                                                                                                                                                                                                                                                                                                                                                                                                                                                                                                                                                                                                                                                                                                                                                                                                                                                                                                                                                                                                                                                                                                                                                                                                                                                                                                                                                                                                                                                                                                                                                                                                                                                                                      
                                                                                                                                                                                                                                                                                                                                                                                                                                                                                                                                                                                                                                                                                                                                                                                                                                                                                                                                                                                                                                                                                                                                                                                                                                                                                                                                                                                                                                                                                                                                                                                                                                                                                                                                                                                                                                                                                                                                                                                                                                                                                                                                      \begin{prop}
Let $\Lambda=\langle a_1\rangle\bot ...\bot\langle a_n\rangle$ be a skew-hermitian lattice such that $N(a_1),...,$ $ N(a_n)\in\pi k^{*2}$, for any prime $\pi$ and $\mu=\nu(4)$. Then $H(\Lambda)=k^*.$                                                                                                                                                                                                                                                                                                                                                                                                                                                                                                                                                                                                                                                                                                                                                                                                                                                                                                                                                                                                                                                                                                                                                                                                                                                                                                                                                                                                                                                                                                                                                                                                                                                                                                                                                                                                                                                                                                                                                                                                                                                                                                                                   
\label{nu4}
\end{prop}                                                                                                                                                                                                                                                                                                                                                                                                                                                                                                                                                                                                                                                                                                                                                                                                                                                                                                                                                                                                                                                                                                                                                                                                                                                                                                                                                                                                                                                                                                                                                                                                                                                                                                                                                                                                                                                                                                                                                                                                                                                                                                                    
                                                                                                                                                                                                                                                                                                                                                                                                                                                                                                                                                                                                                                                                                                                                                                                                                                                                                                                                                                                                                                                                                                                                                                                                                                                                                                                                                                                                                                                                                                                                                                                                                                                                                                                                                                                                                                                                                                                                                                                                                                                                                                                                 \begin{dem}
                                                                                                                                                                                                                                                                                                                                                                                                                                                                                                                                                                                                                                                                                                                                                                                                                                                                                                                                                                                                                                                                                                                                                                                                                                                                                                                                                                                                                                                                                                                                                                                                                                                                                                                                                                                                                                                                                                                                                                                                                                                                                                                                It suffices to prove the result when $n=2$. Hence, we can assume that $\nu(a_2)-\nu(a_1)=\nu(4)$. By Lemma \ref{ret}, we can suppose that the lattice $\langle a_1\rangle\bot\langle a_2\rangle$, has the form $L=\langle q\rangle\bot\langle\epsilon q\rangle$, where $\epsilon \in k^*$ and $\nu(\epsilon)=\nu(4)$. Moreover, we can assume that $q$ is prime and, since all units in $k$ are norms from $k(j)$, we can assume that $q\in ik(j)$. We prove that $H(L)=k^*$ from where the result follows. Because of Corollary \ref{kest2}, it suffices to prove that $L$ represents a prime element whose norm does not belong to the square class of $\pi.$ Remember that if $q\in ik(j)$, then $q\al=\bar{\al}q$ for any $\al\in k(j).$                                                                                                                                                                                                                                                                                                                                                                                                                                                                                                                                                                                                                                                                                                                                                                                                                                                                                                                                                                                                                                                                                                                                                                                                                                                                                                                                                                                                                                                                                                                                                                                                                                                                                                                                                                                                                                                                                                                                                                                                                                                                                                                                                                                                                                                                                                                                                                                                                                                                                                                                                                                                                                                                                                                                                                                                                                                                                                                                                                                                                                                                                                                                                                                                                                                                                                                                                                                                                                                                                                                                                                                                                                                                                                                                                                                                                                                                                                                                                                                                                                                                                                                                                                                                                                                                  We compute, for $\alpha\in \OO_{k(j)}$: $$\left(\begin{array}{cc} 1 & \al \end{array}\right)\left(\begin{array}{cc} q & 0 \\ 0 & \epsilon q \end{array}\right)\left(\begin{array}{cc} 1 \\ \bar{\al} \end{array}\right)=(1+\epsilon \al^2)q,$$ where $N(1+\epsilon \al^2)=1+\epsilon T(\alpha^2)+\epsilon^2N(\al^2)$. Now, since $\nu(\epsilon)=\nu(4)$ and $\epsilon\in k^*$, we have $\epsilon=4\varepsilon,$ with $\varepsilon \in \OO_k^*$. Hence, we conclude $N(1+\epsilon\al^2)=1+4\varepsilon T(\al^2)+16\varepsilon^2N(\al^2).$ Then, in virtue of the Local Square Theorem \cite[\S 63]{O}, it is enough to find $\al$ such that $1+4\varepsilon T(\al^2)$ is not a square. In fact, it is known \cite[\S 63]{O} that there exists a unit of minimal quadratic defect of the form $\Delta=1+4\be\; (\be \in \OO_k^*)$.                                                                                                                                                                                                                                                                                                                                                                                                                                                                                                                                                                                                                                                                                                                                                                                                                                                                                                                                                                                                                                                                                                                                                                                                                                                                                                                                                                                                                                                                                                                                                                                                                                                                                                                                                                                                                                                                                                                                                                                                                                                                                                                                                                                                                                                                                                                                                                                                                                                                                                                                                                                                                                                                                                                                                                                                                                                                                                                                                                                                                                                                                                                                                                                                                                                                                                                                                                                                                                                                                                                                                                                                                                                                                                                                                                                                                                                                                                                                                                                                                                                                                                Hence, it suffices to choose $\al\in \OO_{k(j)}$ such that $T(\alpha^2)=\varepsilon^{-1}\be$. Since $\OO_{k(j)}=\OO_k[\om],$ where $\om=\frac{1+j}{2}$ (\S \ref{ad}), if $\eta\in\OO_{k(j)}$, with $\eta=a+b\om$ for $a,b\in\OO_k$, then $T(\eta)=2a+b$. In particular,  $T(\varepsilon^{-1}\be\om)=\varepsilon^{-1}\be.$ Now, there exists $\al\in \OO_{k(j)}$ such that $\al^2\equiv \varepsilon^{-1}\be\om$ mod($\pi$) since the residue field $\OO_{k(j)}/\pi\OO_{k(j)}$ is perfect of characteristic 2. Hence $T(\al^2)\equiv \varepsilon^{-1}\be$ mod($\pi$) and $1+4\varepsilon T(\al^2)\equiv 1+4\be$ mod($4\pi$). By Local Square Theorem, we have $1+4\varepsilon T(\al^2)= (1+4\be)u^2=\Delta u^2$, for some $u\in \OO_k^*$. Therefore,  $1+4\varepsilon T(\al^2)$ is not a square. We conclude that $N(1+\epsilon\al^2)$ is not a square. Then $N\big(k(i)^*\big)\neq N\Big(k\big((1+\epsilon\al^2)q\big)^*\Big)$. We conclude that $H(\Lambda)=k^*$ as stated.
\end{dem}                                                                                                                                                                                                                                                                                                                                                                                                                                                                                                                                                                                                                                                                                                                                                                                                                                                                                                                                                                                                                                                                                                                                                                                                                                                                                                                                                                                                                                                                                                                                                                                                                                                                                                                                                                                                                                                                                                                                                                                                                                                                                                                                                                                                                                                                                                                                                                                                                                                                                                                                                                                                                                                                                                                                                                                                                                                                                                                                                                                                                                                                                                                                                                                                                                                                                                                                                                                                                                                                                                                                                                                                                                                                                                                                                                                                                                                                                                                                                                                                                                                                                                                                                                                                                                                      

The procedure above cannot be extended to the case $\mu>\nu(4)$ because, in that case, $N(1+\epsilon\al^2)$ is a square. These cases are treated only when the base field is $k=\Q_2$ by the methods used for Case I.
                                                                                                                                                                                                                                                                                                                                                                                                                                                                                                                                                                                                                                                                                                                                                                                                                                                                                                                                                                                                                                                                                                                                                                                                                                                                                                                                                                                                                                                                                                                                                                                                                                                                                                                                                                                                                                                                                                                                                                                                                                                                                                                                      
By the discussion at the beginning of the section, in rank 2 case, we consider lattices $\La$ of the form $\langle i\rangle\bot\langle\ep i\rangle$, where $\ep\in\Q_2^*$ and $\nu(4)<\nu(\ep)\leq\nu(16).$ In this case, the computation depends on the uniformizing parameter. For every prime $\pi$, set $i_{_\pi}\in \OO_{k(j)}i$ such that $i_{_\pi}^2=\pi.$ Remember that, if we prove that $H(\La)=\Q_2^*$, then $H(\La)=\Q_2^*$ for lattices $\La$ of arbitrary rank. Hence, from Lemma \ref{kestpi} and Theorem \ref{eq} next proposition follows.  
                                                                                                                                                                                                                                                                                                                                                                                                                                                                                                                                                                                                                                                                                                                                                                                                                                                                                                                                                                                                                                                                                                                                                                                                                                                                                                                                                                                                                                                                                                                                                                                                                                                                                                                                                                                                                                                                                                                                                                                                                                                                                                                       \begin{prop}
Let $\Lambda=\langle a_1\rangle\bot ...\bot\langle a_n\rangle$ be a skew-hermitian lattice such that $N(a_1)$,..., $N(a_n)\in \pi\Q_2^{*2}$ and  $\mu$ satisfies $\nu(4)<\mu\leq \nu(16)$. Then $H(\Lambda)=\Q_2^*.$                                                                                                                                                                                                                                                                                                                                                                                                                                                                                                                                                                                                                                                                                                                                                                                                                                                                                                                                                                                                                                                                                                                                                                                                                                                                                                                                                                                                                                                                                                                                                                                                                                                                                                                                                                                                                                                                                                                                                                                                                                                                                                                                   
\label{propcaso2}
\end{prop}                                                                                                                                                                                                                                                                                                                                                                                                                                                                                                                                                                                                                                                                                                                                                                                                                                                                                                                                                                                                                                                                                                                                                                                                                                                                                                                                                                                                                                                                                                                                                                                                                                                                                                                                                                                                                                                                                                                                                                                                                                                                                                                    
\section{Examples}
 
\begin{enumerate}

\item Let $\La=\langle i+j\rangle\bot\langle 8(i+j)\rangle$ be a skew-hermitian lattice, where $D=\left(\frac{2,5}{\Q}\right)$, i.e., $i^2=2$ and $j^2=5$. $D$ ramifies only at 2 and 5. Moreover, $\La$ is unimodular for $p\neq 2,7$. Hence, $H(\La_p)=\Z_p^*\Q_p^{*2}$ for $p\neq 2,5,7$ \cite{K1} and $H(\La_5)=\Z_5^*\Q_5^{*2}$ \cite{B}. Then, the spinor class field $\Sigma_\La$ satisfies $\Sigma_\La\subset\Q(\sqrt{-1},\sqrt{2},\sqrt{7})$. On the other hand, since the algebra is decomposed at infinite and the associated quadratic form is indefinite, we have that the class and the spinor genus coincide and $\Sigma_\La$ satisfies $\Sigma_\La\subset \R$. From here and the fact that $H(\La_2)=N\big(\Q_2(i+j)^*\big)=N\big(\Q_2(\sqrt{7})^*\big)$ (see Table \ref{completa}), we conclude that $\Sigma_\La\subset\Q(\sqrt{7})$. For $p=7$ the algebra decomposes and the quadratic form associated to $\langle i+j\rangle$ has discriminant $N(i+j)=-7$. Then, the quadratic lattice associated to $\La$ has a decomposition of the type $$\langle a\rangle\bot\langle -7a\rangle\bot\langle 8a\rangle\bot\langle -56a\rangle,$$ where $a\in\Z_7^*$. Hence, if $\langle a\rangle\bot\langle 8a\rangle=\Z_7 e_1\bot\Z_7 e_2$, we consider the rotation $\tau_{e_1}\tau_{3e_1+e_2}$, for which $\theta(\tau_{e_1}\tau_{3e_1+e_2})=3\Q_7^{*2}.$ This fact tell us that $\Sigma_\La$ is decomposed at 7 because $3\notin N\big(\Q_7(\sqrt{7})^*\big).$ We conclude that $\Sigma_\La=\Q$ and therefore, the class number of $\La$ is 1. 
 
 
\item Let consider the family of lattices $\La=\langle i\rangle\bot\langle 2^ti\rangle$, for $t>0$, where $D=\left(\frac{2,5}{\Q}\right)$. $D$ ramifies only at 2 and 5. The lattice $\La$ is unimodular for $p\neq 2$. We have that $H(\La_p)=\Z_p\Q_p^{*2}$ for $p\neq 2$, in virtue of the computations in \cite{K1} (for $p\neq 5$) and \cite[Theorem 4]{B} (for $p=5$). Hence, the spinor class field $\Sigma_\La$ can ramify only at 2 and $\infty$. So, $\Sigma_\La\subset\Q(\sqrt{-1},\sqrt{2})$. Observe that the algebra $D$ decomposes at infinity and the quadratic form corresponding to $\La$ is indefinite. Hence, class and spinor genus of $\La$ coincide and $\Sigma_\La\subset\R$. On the other hand, for $p=2$, Table \ref{completa} tell us that $H(\La_2)=\Q_2^*$ if $t\leq 4$ and $H(\La_2)=N\big(\Q_2(i)^*\big)$ if $t>4$, whence $\Sigma_\La$ decomposes at 2 for $t\leq 4$ and ramifies at 2 for $t>4$. We conclude that $\Sigma_\La=\Q$ for $t\leq 4$, while $\Sigma_\La=\Q(\sqrt{2})$ for $t>4$. In the first case, Hasse principle holds for $\La$. In the second case, the class number of $\La$ is 2.  
\end{enumerate}

\section*{Acknowledgement}

The research was partly supported by Fondecyt, Project No. 1120565.

\bibliographystyle{amsplain}
\bibliography{spanish}







\end{document}